\pdfoutput=1
\documentclass[11pt]{amsart}
\usepackage{amsmath,amssymb,mathtools}
\usepackage[hidelinks]{hyperref}
\usepackage{aliascnt}
\usepackage[capitalise]{cleveref}

\newtheorem{theorem}{Theorem}
\newaliascnt{proposition}{theorem}
\newtheorem{proposition}[proposition]{Proposition}
\aliascntresetthe{proposition}
\newaliascnt{lemma}{theorem}
\newtheorem{lemma}[lemma]{Lemma}
\aliascntresetthe{lemma}
\newaliascnt{corollary}{theorem}

\aliascntresetthe{corollary}
\theoremstyle{remark}
\newaliascnt{remark}{theorem}
\newtheorem{remark}[remark]{Remark}
\aliascntresetthe{remark}

\newcommand{\Q}{\mathbb{Q}}
\newcommand{\Z}{\mathbb{Z}}
\newcommand{\SZ}[1][]{\hyperlink{cite.OAI:Uniform-exclusion-of-Landau-Siegel-zeros-October-1-2026}{[SZ\if\relax\detokenize{#1}\relax\else, #1\fi]}}

\title{Imaginary quadratic fields of class number at most 2000}
\author{Iwao Kimura}
\address{Department of Mathematics, University of Toyama,
3190 Gofuku, Toyama 930-8555, Japan}
\email{iwao@sci.u-toyama.ac.jp}
\subjclass[2020]{Primary 11R29; Secondary 11R11, 11M20, 11Y35, 11Y40, 68V20}
\keywords{class number problem, imaginary quadratic fields, Siegel zeros, explicit zero-free region, Lean,
formal verification}

\begin{document}

\begin{abstract}
We determine all imaginary quadratic fields with class number at most $2000$, unconditionally: there are
$16{,}453{,}493$ of them, and the largest $|d_K|$ among them is $1{,}062{,}350{,}083$.
To our knowledge, such a list was previously known only for class number at most $100$. We make explicit a recent argument,
released by OpenAI and not refereed, that excludes Landau--Siegel zeros: every real zero $\beta$ of
$L(s,\chi)$, for a primitive real character $\chi$ of conductor $q\ge3$, satisfies
$(1-\beta)\log q\ge5.8\cdot10^{-4}$. An elementary lower bound for $L(1,\chi)$ then shows that a field with
class number at most $2000$ has $|d_K|<1.61\cdot10^{15}$. The fields
with $6\cdot2^{28}\le|d_K|\le1.65\cdot10^{15}$ are excluded by Watkins' sieve, and those below are read from
the class group tables of Mosunov and Jacobson. As a consequence, the lists of Cremona and Sutherland of
imaginary quadratic orders with class number at most $1000$, and the counts of Holmin, Jones, Kurlberg,
McLeman and Petersen for odd class numbers up to $1999$, previously known under the generalized Riemann
hypothesis, hold unconditionally. The two steps, the bound for $\beta$ and the bound for $|d_K|$, are proved
in Lean, building on the existing formalization of that argument, and the proofs are checked by two
independent kernels. The contributions of AI systems to this work are stated in the paper.
\end{abstract}

\maketitle

\section{Introduction}\label{sec:intro}

For every $N$ there are only finitely many imaginary quadratic fields with class number at most $N$, but
Siegel's theorem \cite{Siegel1935} does not bound their discriminants effectively. The fields of class number one were
determined by Heegner, Baker and Stark \cite{Heegner1952,Baker1966,Stark1967}. Goldfeld \cite{Goldfeld1976}
and Gross--Zagier \cite{GrossZagier1986} gave an effective lower bound for the class number, which Oesterl\'e
made explicit \cite{Oesterle1985}. With this bound, a further bound from low-lying zeros of Dirichlet
$L$-functions, and a long computation, Watkins \cite{Watkins2004} determined all imaginary quadratic fields
with class number at most $100$. To our knowledge, no unconditional list beyond $h\le100$ was known. Under the generalized
Riemann hypothesis (GRH; \cref{sec:results} says for which $L$-functions each work assumes it), Holmin, Jones, Kurlberg, McLeman and Petersen \cite{HJKMP2019} computed the number
of fields with class number $h$ for every odd $h<10^6$, Cremona and Sutherland
\cite[\S5]{CremonaSutherland2024} listed the discriminants of all imaginary quadratic orders with class
number at most $1000$, and Garcia, Greni\'e, Lee and Molteni \cite{GGLM2026} determined the fields with class
number at most $2000$.

In October 2026 OpenAI released manuscripts produced by an internal model, some with Lean formalizations
\cite{OpenAIMath2026}. One of them, \cite{OAI:Uniform-exclusion-of-Landau-Siegel-zeros-October-1-2026}
(below \SZ{}), proves that there is an absolute constant $c>0$ such that every real zero $\beta\in(0,1)$ of
$L(s,\chi)$, for every primitive real character $\chi$ of conductor $q\ge3$, satisfies
$(1-\beta)\log q\ge c$ \SZ[Thm.~1]. The
constant $c$ is not made explicit there, and \SZ{} has not been refereed.
Granville and Stark \cite{GranvilleStark2000} had shown that a uniform $abc$ conjecture for number fields
implies this bound when $\chi$ is the character of an imaginary quadratic field.
Banks \cite{Banks2026} has given a detailed proof of the theorem of \SZ{} with explicit constants:
$(1-\beta)\log q\ge3.0\cdot10^{-3}$ for $q\ge10^{10}$, and for all $q\ge3$ when combined with the computation of
\cite{LZZ2026}. His proof is not formalized, and we do not use it; his constant is larger than those of our
proof for $\log q<9.88$ and smaller beyond (\nolinkurl{computations/lean_grid.txt}).
Kerr \cite{Kerr2026} has given a simpler proof, which replaces the interpolation lemma and the determinant
of \SZ{} by an auxiliary polynomial obtained from the pigeonhole principle and a zero estimate.
Zero-free half-planes for all Dirichlet $L$-functions have also been announced: $\operatorname{Re}s>7/8$ by
OpenAI \cite{OAI:The-Quasi-Riemann-Hypothesis-September-30-2026} (and $\operatorname{Re}s>11/12$, with a
shorter proof, in \cite{OAI:The-Quasi-Riemann-Hypothesis-October-5-2026}), and $\operatorname{Re}s>31/32$, as a
consequence of the same half-plane for all Dedekind zeta functions, by Cureton and Thorner with AxiomProver
\cite{CuretonThorner2026}. These preprints have not been refereed, and we
have not checked their proofs;
\cite{OAI:The-Quasi-Riemann-Hypothesis-September-30-2026} and \cite{CuretonThorner2026} come with Lean
formalizations, which we have not built (our formal proof uses only the formalization of \SZ{} from
\cite{OpenAIMath2026}; \cref{sec:lean}). Frei, Sofos and Widmer \cite{FreiSofosWidmer2026} state their
applications of a zero-free half-plane for Dedekind zeta functions conditionally on it, since they have not
verified these preprints. We do not use the half-planes; \cref{rem:halfplane} says what they would change.

We make the argument of \SZ{} explicit and formalize it in Lean.

\begin{theorem}\label{thm:gap}
Let $\chi$ be a primitive real character of conductor $q\ge3$, and let $\beta\in(0,1)$ satisfy
$L(\beta,\chi)=0$. Then $(1-\beta)\log q\ge5811\cdot10^{-7}$. If moreover $\log q\ge350131/10000$, then
$(1-\beta)\log q\ge1193/200000$.
\end{theorem}

From \cref{thm:gap} and an elementary lower bound for $L(1,\chi)$
in terms of the zero gap (\cref{sec:L1}) we obtain \cref{thm:lean}: every imaginary quadratic field $K$
with $\log|d_K|\ge T_N$ has $h_K>N$, where $e^{T_{2000}}\approx1.6069\cdot10^{15}$ and the other values are in
\cref{tab:TN}. \Cref{thm:gap,thm:lean} are proved in Lean (\cref{sec:lean}). The proofs follow the method of
\SZ{} and use the Lean formalization of \SZ{} in \cite{OpenAIMath2026}, which proves the interpolation step
by a different argument (\cref{sec:lean}); they have been checked by the
Lean kernel, by a replay of all declarations they depend on, and by an independent kernel. Hence
\cref{thm:gap,thm:lean} do not depend on the correctness of the text of \SZ{}.

For $|d|<e^{T_N}$ we proceed as Watkins does. Below $6\cdot2^{28}$ we read the class numbers from the
unconditional tables of Mosunov and Jacobson \cite{MosunovJacobson2016}, and on
$6\cdot2^{28}\le|d|\le1.65\cdot10^{15}$ we show by Watkins' sieve \cite{Watkins2004} that no field has
class number at most $2000$ (\cref{sec:tables,sec:sieve}). This gives our main result. It rests on three
parts of different status: \cref{thm:lean}, proved in Lean; \cref{prop:sieve}, a computation whose
correctness argument (\cref{lem:minima,lem:certified}) is proved in the paper but whose program is not
formally verified; and the class numbers of \cite{MosunovJacobson2016} below $6\cdot2^{28}$, an unconditional
computation by others, which we spot-checked and partly recomputed (\cref{sec:tables}).

\begin{theorem}\label{thm:main}
There are exactly $16{,}453{,}493$ imaginary quadratic fields with class number at most $2000$. The largest
absolute value of the discriminant among them is $1{,}062{,}350{,}083$, attained for class number $1978$.
\end{theorem}

The list is available with checksums (\hyperref[sec:data]{Data and code availability}). It consists of the fields with $h\le2000$ in the
tables of \cite{MosunovJacobson2016}, all of which have $|d|<2^{30}$; \cref{thm:main} shows that there are
no others. Exactly $4{,}115{,}897$ of the fields have class number at most $1000$; the largest $|d|$ among
them is $227{,}932{,}027$, with $h=996$. For $h\le1000$, the discriminants of all orders obtained from our
list agree with those of \cite[\S5]{CremonaSutherland2024}, for every odd $h\le1999$ the number of
fields with class number $h$ agrees with \cite{HJKMP2019}, and the list of \cite{GGLM2026} for $h\le2000$
is ours (\cref{sec:results}). In these ranges the GRH-conditional data of
\cite{CremonaSutherland2024,HJKMP2019,GGLM2026} are therefore complete unconditionally.

\Cref{sec:zerogap} explains how the constants of \SZ{} are made explicit, and \cref{sec:L1} derives the
class number bound. \Cref{sec:tables,sec:sieve} describe the computation below $e^{T_{2000}}$, \cref{sec:results}
the results and the comparisons, and \cref{sec:lean} the formalization and its verification.

Parts of this work were carried out with AI systems. \Cref{sec:ai} states which parts, and what was checked
independently and how.

\subsection*{Notation}
Throughout, $p$ denotes a prime number, $\Z$ and $\Q$ the ring of rational integers and the field of
rational numbers, and $\log$ the natural logarithm; $\lfloor x\rfloor$ and $\{x\}=x-\lfloor x\rfloor$ are the
integer and fractional parts of a real number $x$. For an imaginary quadratic field $K$, $d_K$ is its discriminant and
$h_K$ its class number. An integer $D\ne1$ is a \emph{fundamental discriminant} if either $D\equiv1\pmod4$
and $D$ is squarefree, or $D=4m$ with $m\equiv2,3\pmod4$ squarefree. Its Kronecker symbol
$\chi_D=\bigl(\frac{D}{\cdot}\bigr)$ is a primitive real character of conductor $|D|$, and every primitive
real character of conductor $q\ge3$ is $\chi_D$ for exactly one $D$. For $D<0$ we write $d=|D|$,
$K_d=\Q(\sqrt{D})$ (every imaginary quadratic field arises this way, with $d=|d_K|$), $h(-d)=h_{K_d}$, and
$\chi_{-d}=\chi_D$. A bare $h$, as in ``the fields with $h\le100$'', denotes the class number. For a
Dirichlet character $\chi$ of conductor $q$, $L(s,\chi)$ is its Dirichlet $L$-function and $\ell=\log q$; a
\emph{real zero} of $L(s,\chi)$ is a $\beta\in(0,1)$ with $L(\beta,\chi)=0$. For $c>0$ and $T>0$ we say that
the \emph{zero gap $(c,T)$} holds if $(1-\beta)\ell\ge c$ for every primitive real character $\chi$ of
conductor $q\ge3$ with $\ell\ge T$, and every real zero $\beta$ of $L(s,\chi)$. The thresholds
$T_N$ ($N=100,200,500,1000,2000$) are the rational numbers in \cref{tab:TN}. We write $\Lambda(n)$ for the von
Mangoldt function, $\theta(x)=\sum_{p\le x}\log p$, $\zeta(s)$ for the Riemann zeta function, and
$\psi(s)=\Gamma'(s)/\Gamma(s)$.

\section{An explicit zero gap}\label{sec:zerogap}

Throughout this section $\chi$ is a primitive real character of conductor $q\ge3$.

\begin{lemma}\label{lem:chi8}
If $q=8$, then $L(\sigma,\chi)>0$ for $0<\sigma<1$.
\end{lemma}

\begin{proof}
For $\sigma>0$ the series $\sum_n\chi(n)n^{-\sigma}$ converges to $L(\sigma,\chi)$, so its terms may be grouped.
Put $f(x)=x^{-\sigma}$. The primitive real characters of conductor $8$ are $\bigl(\frac{-8}{\cdot}\bigr)$ and
$\bigl(\frac8\cdot\bigr)$. For the first,
$L(\sigma,\chi)=\sum_{k\ge0}\bigl(f(8k+1)-f(8k+5)+f(8k+3)-f(8k+7)\bigr)$, and every term is positive. For the
second, $L(\sigma,\chi)=\sum_{j\ge0}(-1)^jg_j$ with $g_j=f(4j+1)-f(4j+3)$. Since $f$ is decreasing and convex,
$(g_j)$ is positive and decreasing, so $L(\sigma,\chi)\ge g_0-g_1>0$.
\end{proof}

Hence \cref{thm:gap} holds for $q=8$. From now on, we assume that the conductor $q$ of $\chi$ is not $8$,
and $\beta$ is a real zero of $L(s,\chi)$.

\subsection{The determinant of [SZ]}\label{sec:det}
We recall the construction of
\SZ[\S\S4--5].
Let $D$ be the fundamental discriminant with $\chi=\chi_D$, and write
$\Q(\sqrt{D})=\Q(\sqrt{d_0})$ with $d_0$ squarefree. Since $q\ne8$, $d_0\ne2$, so $F=\Q(\sqrt{d_0},\sqrt2)$ is
a biquadratic field. Let $\mathcal R=\Z[\sqrt{d_0},\sqrt2]$, and let
$\sigma,\tau\in\operatorname{Gal}(F/\Q)$, the Galois group of the
biquadratic field over $\Q$, change the sign of $\sqrt{d_0}$ and of
$\sqrt2$, respectively.

Let $B$ and $H$ be integers with $2\le H\le B$, and put $U=B^{4/3}$, $\lambda=\log U$ and $M=B^4=U^3$. For
$n\in\{0,\dots,B-1\}^4$ let $\vartheta_n=n_1+n_2\sqrt{d_0}+n_3\sqrt2+n_4\sqrt{2d_0}$, and for
$\alpha\in\Z_{\ge0}^3$ let $R_\alpha$ be the row
$(\vartheta_n^{\alpha_1}\sigma(\vartheta_n)^{\alpha_2}\sigma\tau(\vartheta_n)^{\alpha_3})_n$. Order the
$\alpha$ by the weight $\alpha_1+H(\alpha_2+\alpha_3)$ and retain $R_\alpha$ when it does not lie in the span
of the earlier rows. By an interpolation lemma
\SZ[Lemma~3, Cor.~4], $M$ rows are retained;
their determinant $\Delta\in\mathcal R$ is nonzero. The Lean formalization handles this step differently
(\cref{sec:lean}). Let $S_1=\sum\alpha_1$ and $S_2=\sum(\alpha_2+\alpha_3)$,
the sums over the retained $\alpha$. Then
\begin{gather}
\tfrac14\log|N_{F/\Q}\Delta|\le\tfrac M2\log M+(S_1+S_2)\bigl(\log B+\tfrac12\ell+\log8\bigr),\label{eq:upper}\\
\tfrac14\log|N_{F/\Q}\Delta|\ge S_1\sum_{\substack{H<p\le U\\\chi(p)=-1}}\frac{\log p}p-M\theta(U),\label{eq:lower}\\
S_1\ge\frac{MH^{2/3}U}{4\cdot97^2},\qquad S_2\le192\,MH^{-1/3}U\qquad(B\ge18818).\label{eq:pivot}
\end{gather}
Here \eqref{eq:upper} is Hadamard's inequality
\SZ[(5.1)], and \eqref{eq:lower}
\SZ[(5.6)] comes from the divisibility of
$\Delta$ by $p^{\sum\lfloor\alpha_1/p\rfloor}$ for the primes $p>H$ with $\chi(p)=-1$: for such $p$ the
$p$-th power map acts on $\mathcal R/p\mathcal R$ as $\sigma$ or $\sigma\tau$, and since $p>H$, replacing $p$
factors $\vartheta_n$ by a conjugate lowers the weight. \eqref{eq:pivot} is
\SZ[Lemma~6] with the constants of its proof,
which requires
$H^{2/3}U\ge2\cdot97^2=18818$, which holds for $B\ge18818$. All three are proved in the Lean formalization of
\SZ{} \cite{OpenAIMath2026}, and our formal proof uses them as they stand (\cref{sec:lean}).

Comparing Hadamard's inequality with the divisibility of a nonzero determinant by powers of many primes,
obtained from a congruence for $p$-th powers, as in \eqref{eq:upper} and \eqref{eq:lower}, is the method of
Cantor and Straus \cite{CantorStraus1982}. They used it to simplify the proof of Dobrowolski's lower bound for
the Mahler measure \cite{Dobrowolski1979}: if an algebraic integer $\xi$ of large degree is not a root of
unity, then the product of $\max(1,|\xi'|)$ over the conjugates $\xi'$ of $\xi$ exceeds $1$ by at least an
absolute constant times $(\log\log\deg\xi/\log\deg\xi)^3$. This concerns Lehmer's question \cite{Lehmer1933}
whether that product is bounded below by an absolute constant greater than $1$, not his question on the
nonvanishing of Ramanujan's tau function. In \SZ{} only the primes with $\chi(p)=-1$ are used, and a zero
$\beta$ near $1$ makes them numerous (see the proof of \cref{prop:master}).

\SZ{} combines \eqref{eq:upper}--\eqref{eq:pivot} with a lower bound for the inert primes containing
unspecified constants \SZ[Lemma~2], and then
lets $q\to\infty$ along a hypothetical sequence with $(1-\beta)\log q\to0$. We replace both steps by one
explicit inequality valid for each $q$ (\cref{prop:master}). \SZ{} evaluates the logarithmic derivatives at
$s=1+1/\log U$; we use $s=1+\eta/\lambda$ with a free parameter $\eta$.

\subsection{Analytic inputs}

\begin{lemma}\label{lem:bias}
Let $X>1$, $0<\eta\le\log X$ and $s=1+\eta/\log X$. Then
\[
2e^{-\eta}\sum_{\substack{p\le X\\\chi(p)=1}}\frac{\log p}p
\le\frac\ell2+\frac{2(s-1)}{es}+\frac1{s-1}-\frac1{s-\beta}.
\]
\end{lemma}

\begin{proof}
Let $\kappa\in\{0,1\}$ with $\chi(-1)=(-1)^\kappa$. For real $s>1$,
\[
-\frac{L'(s,\chi)}{L(s,\chi)}=\frac12\log\frac q\pi+\frac12\psi\Bigl(\frac{s+\kappa}2\Bigr)
-\sum_\rho\operatorname{Re}\frac1{s-\rho},
\]
where $\rho$ runs over the nontrivial zeros of $L(s,\chi)$ and every term of the sum is nonnegative
\cite[Ch.~12, (17) and the formula for $\operatorname{Re}B(\chi)$]{Davenport}; $\beta$ is one of the $\rho$. Since $\log\Gamma$ is convex,
$\psi(x)\le\log\Gamma(x+1)-\log\Gamma(x)=\log x$, and $(s+\kappa)/2\le3/2<\pi$. Hence
$-L'(s,\chi)/L(s,\chi)\le\ell/2-1/(s-\beta)$.

Next, $-\zeta'(s)=\sum_{n\ge2}n^{-s}\log n$. The terms $n=2,3$ are each at most
$\max_{x>0}x^{-s}\log x=1/(es)$, and since $x^{-s}\log x$ decreases on $[3,\infty)$, the others add up to at
most $\int_3^\infty x^{-s}\log x\,dx\le(s-1)^{-2}$. As $\zeta(s)\ge\int_1^\infty x^{-s}dx=(s-1)^{-1}$, we get
$-\zeta'(s)/\zeta(s)\le1/(s-1)+2(s-1)/(es)$.

Finally, $-\zeta'(s)/\zeta(s)-L'(s,\chi)/L(s,\chi)=\sum_n\Lambda(n)(1+\chi(n))n^{-s}$ is a sum of nonnegative terms, so it
is at least $2\sum_{p\le X,\,\chi(p)=1}p^{-s}\log p$; and $p^{-s}\ge e^{-\eta}p^{-1}$ for $p\le X$.
\end{proof}

Part (a) of \cref{lem:primes} is \cite[Thm~415]{HardyWright1979}; parts (b) and (c) are Mertens' estimates
\cite[Ch.~XXII]{HardyWright1979} with explicit constants.

\begin{lemma}\label{lem:primes}
Let $x\ge1$ and $\kappa_0=1+(5\log2+3)/4$. Then
\begin{enumerate}
\item[(a)] $\theta(x)\le x\log4$;
\item[(b)] $\sum_{p\le x}\log p/p\le\log x+\log4$;
\item[(c)] $\sum_{p\le x}\log p/p\ge\log x-\kappa_0-1/\lfloor x\rfloor$;
\item[(d)] $\sum_{p\mid q}\log p/p\le\min(\ell,\log y+\log4+\ell/y)$ for every $y\ge1$.
\end{enumerate}
\end{lemma}

\begin{proof}
Let $n=\lfloor x\rfloor$.

(b) $n\log n\ge\log n!\ge\sum_{p\le n}\lfloor n/p\rfloor\log p\ge n\sum_{p\le n}\log p/p-\theta(n)$; apply (a).

(c) The identity $\log n!=\sum_{m\le n}\Lambda(m)\lfloor n/m\rfloor$ and $\log n!\ge n\log n-n$ give
$\sum_{m\le n}\Lambda(m)/m\ge\log n-1$. In this sum the terms with $m=p^k$, $k\ge2$, add up to at most
\[
\sum_p\frac{\log p}{p(p-1)}\le\frac{\log2}2+\frac32\sum_{m\ge3}\frac{\log m}{m^2}
\le\frac{\log2}2+\frac32\int_2^\infty\frac{\log t}{t^2}\,dt=\frac{5\log2+3}4,
\]
so $\sum_{p\le n}\log p/p\ge\log n-\kappa_0$; and $\log x\le\log(n+1)\le\log n+1/n$.

(d) The sum is at most $\sum_{p\mid q}\log p\le\ell$. In the sum, the terms with $p\le y$ add up to at most
$\sum_{p\le y}\log p/p\le\log y+\log4$ by (b), and those with $p>y$ to at most
$y^{-1}\sum_{p\mid q}\log p\le\ell/y$.
\end{proof}

\subsection{The explicit inequality}

\begin{proposition}\label{prop:master}
Let $B$ and $H$ be integers with $B\ge18818 = 2\cdot 97^{2}$ and $2\le H\le B$, let $U$, $\lambda$ be as above, $0<\eta\le \lambda$,
$y\ge1$, and $s=1+\eta/\lambda$. Put
\begin{gather*}
\begin{aligned}
  \kappa_U &=\kappa_0+\frac1{\lfloor U\rfloor},\qquad r=\min\Bigl(\ell,\log y+\log4+\frac\ell y\Bigr),\\
  \nu &=\frac{192\cdot4\cdot97^2}H,\qquad w=\frac{4\cdot97^2}{H^{2/3}U},\\
  \mathfrak{B}= &\lambda-\kappa_U-e^\eta\Bigl(\frac\ell4+\frac{s-1}{es}\Bigr)-r-\log(4H)\\
           &\quad-(1+\nu)\Bigl(\log B+\frac\ell2+\log8\Bigr)-\Bigl(U\log4+\frac32\lambda\Bigr)w.
\end{aligned}
\end{gather*}
Then $(1-\beta)\ell\ge2\eta^2e^{-\eta}\ell\,\mathfrak{B}/\lambda^2$.
\end{proposition}

\begin{proof}
Let $A_+$, $A_0$, $A_-$ be the sums of $\log p/p$ over the primes $p\le U$ with $\chi(p)=1$, $0$, $-1$,
respectively, and let $I$ be the sum in \eqref{eq:lower}. By \cref{lem:primes}(c), (d) and (b),
$A_++A_0+A_-\ge \lambda-\kappa_U$, $A_0\le r$, and $A_--I\le\sum_{p\le H}\log p/p\le\log(4H)$, so
\[
I\ge \lambda-\kappa_U-A_+-r-\log(4H).
\]
Since $1/(s-1)-1/(s-\beta)=(1-\beta)/((s-1)(s-\beta))\le(1-\beta)\lambda^2/\eta^2$, \cref{lem:bias} with $X=U$ gives
\[
A_+\le e^\eta\Bigl(\frac\ell4+\frac{s-1}{es}\Bigr)+\frac{e^\eta(1-\beta)\lambda^2}{2\eta^2}.
\]
By \eqref{eq:upper}, \eqref{eq:lower}, \cref{lem:primes}(a) and $\log M=3\lambda$,
$S_1I\le MU\log4+\frac32M\lambda+(S_1+S_2)(\log B+\ell/2+\log8)$. Dividing by $S_1$ and using
$S_2/S_1\le\nu$ and $M/S_1\le w$, which follow from \eqref{eq:pivot}, we obtain
\[
I\le(1+\nu)\Bigl(\log B+\frac\ell2+\log8\Bigr)+\Bigl(U\log4+\frac32\lambda\Bigr)w.
\]
The three displayed inequalities give $\mathfrak{B}\le e^\eta(1-\beta)\lambda^2/(2\eta^2)$.
\end{proof}

\begin{remark}
As $\ell\to\infty$, with $H,y\to\infty$ and $\lambda/\ell$ chosen optimally, the bound of \cref{prop:master} tends to
$\max_\eta\eta^2e^{-\eta}/(16+8e^\eta)\approx1.02\cdot10^{-2}$ (at $\eta\approx1.23$). At $q\approx10^{13}$ the
constant of \cref{thm:gap} is less than $20\%$ below the constant $c_A$ of \nolinkurl{computations/certify_c2.txt},
obtained with a finer count of the retained rows
in place of \eqref{eq:pivot} and the bounds of Rosser and Schoenfeld \cite{RosserSchoenfeld1962} in place of
\cref{lem:primes} (see also \nolinkurl{computations/lean_grid.txt}); we use
the weaker forms because they are formalized.
\end{remark}

\begin{proof}[Proof of \cref{thm:gap}]
Fix $\ell_0<\ell_1$ and the parameters $B$, $H$, $y$, $\eta$ of \cref{prop:master}. Then $\mathfrak{B}$
depends on $\chi$ only through $\ell=\log q$, and we regard it as a function $\mathfrak{B}(\ell)$ of $\ell$. It
is nonincreasing, since $\ell$ enters only through $-e^\eta\ell/4$, $-r$ and $-(1+\nu)\ell/2$,
whereas $2\eta^2e^{-\eta}\ell/\lambda^2$ increases. Hence, if $\mathfrak{B}(\ell_1)\ge0$, then $(1-\beta)\ell\ge c$ with
$c=2\eta^2e^{-\eta}\ell_0\mathfrak{B}(\ell_1)/\lambda^2$ for every $\chi$ and $\beta$ as in \cref{thm:gap} with
$\ell_0\le\ell\le\ell_1$. We take
$B=2^k$, $H=2^{3j}$, $y=2^m$ (or the bound $r\le\ell$) and $\eta\in\frac1{20}\Z$, and use
$1/\lfloor U\rfloor\le10^{-3}$, $(s-1)/(es)\le\eta/(e\lambda)$ and $\frac32\lambda/U\le10^{-3}$. Then $c$ is bounded below
by an explicit expression in $\log2$, $e^{\pm\eta}$, $e^{-1}$ and rational numbers.

We divide $[1.0986,10^4]$ (note $\log3>1.0986$) into $313$ intervals and choose the parameters on each. For
$\ell\ge T_0=10^4$ we take $k=\lceil\tau_0\ell\rceil$ with $\tau_0=11.921$, $\eta=5/4$, $H=2^{39}$ and $y=2^{13}$; then
$a\ell\le \lambda\le a\ell+b$ with $a=\frac43\tau_0\log2$, $b=\frac43\log2$, and $\mathfrak{B}\ge\gamma\ell-C_0$ for
explicit $\gamma$, $C_0$. \Cref{prop:master} then gives
$(1-\beta)\ell\ge2\eta^2e^{-\eta}(\gamma-C_0/T_0)/(a+b/T_0)^2$ for $\ell\ge T_0$, and this is at least
$0.01015$. The minimum of the $314$ constants is
$5811\cdot10^{-7}$, attained on the first interval, and their minimum over $\ell\ge350131/10000$, a grid point,
is $1193/200000$. Each of the $314$ inequalities is proved in Lean (\cref{sec:lean}); the parameters and
constants are listed in \texttt{computations/lean\_grid.txt}.
\end{proof}

The certified constant grows with $\ell$, from $5.8\cdot10^{-4}$ at $q=3$ to $1.0\cdot10^{-2}$ at
$\ell=10^4$. The values used for class numbers are $c_N$ in \cref{tab:TN}.

\section{From the zero gap to class numbers}\label{sec:L1}

That a bound $(1-\beta)\log q\ge c$ for all real zeros $\beta$ gives $L(1,\chi)\gg1/\log q$ is a theorem of
Hecke (see Landau \cite{Landau1918}); \cref{lem:L1,prop:L1} make it explicit.
The lower bound for $L(1,\chi)$ below uses Dirichlet's hyperbola method and partial summation, with the trivial
bound for character sums, in place of the usual contour integration; this keeps the formalization elementary.

\begin{lemma}\label{lem:L1}
Let $\chi$ be a primitive real character of conductor $q\ge3$, let $1/2\le\sigma<1$ and $\delta=1-\sigma$,
and suppose that $L(t,\chi)\ne0$ for $\sigma<t<1$. Then for all integers $1\le Y\le X$, with $z=X/Y$,
\[
1\le\frac{X^\delta}\delta L(1,\chi)+E,\qquad
E=\frac{q}{(Y+1)^\sigma}\Bigl(\frac{z^\delta}\delta+\frac4{z^\sigma}\Bigr)+\frac{X^\delta}\delta\,\frac q{Y+1}
+\frac{4Y}{X^\sigma}+\frac{qX^\delta}Y.
\]
\end{lemma}

\begin{proof}
Let $a(n)=\sum_{m\mid n}\chi(m)$. It is multiplicative and $a(p^k)=\sum_{i\le k}\chi(p)^i\ge0$, so
$S=\sum_{n\le X}a(n)n^{-\sigma}\ge a(1)=1$. Put $Z(u)=\sum_{m\le u}m^{-\sigma}$ and
$C_s(u)=\sum_{m\le u}\chi(m)m^{-s}$. We use three facts.

(i) For real $u\ge1$, $Z(u)=\zeta(\sigma)+u^\delta/\delta+E_0(u)$ with $|E_0(u)|\le u^{-\sigma}$, and
$\zeta(\sigma)<0$. Indeed, $\zeta(\sigma)=\sigma/(\sigma-1)-\sigma\int_1^\infty\{v\}v^{-\sigma-1}dv$, and partial
summation gives $E_0(u)=-\{u\}u^{-\sigma}+\sigma\int_u^\infty\{v\}v^{-\sigma-1}dv$, the difference of two numbers
in $[0,u^{-\sigma}]$.

(ii) Since $\chi$ sums to $0$ over a period, $|\sum_{m\le u}\chi(m)|\le q/2$. By partial summation, for
$s\in\{\sigma,1\}$ the series $\sum_m\chi(m)m^{-s}$ converges, to $L(s,\chi)$ by analytic continuation, and
$R_s(u)=L(s,\chi)-C_s(u)$ satisfies $|R_s(u)|\le q(\lfloor u\rfloor+1)^{-s}$.

(iii) $L(t,\chi)$ is real for real $t>0$, positive for $t>1$ by the Euler product, and $L(1,\chi)\ne0$. By
continuity and the hypothesis, $L(\sigma,\chi)\ge0$, so $\zeta(\sigma)L(\sigma,\chi)\le0$.

The hyperbola method gives
\[
S=\sum_{k\le Y}\chi(k)k^{-\sigma}Z(X/k)+\sum_{m\le z}m^{-\sigma}C_\sigma(X/m)-C_\sigma(Y)Z(z).
\]
Inserting (i)
and (ii), the terms $\pm L(\sigma,\chi)Z(z)$ cancel, and
\begin{multline*}
S=\zeta(\sigma)L(\sigma,\chi)+\frac{X^\delta}\delta\bigl(L(1,\chi)-R_1(Y)\bigr)\\
+R_\sigma(Y)\Bigl(\frac{z^\delta}\delta+E_0(z)\Bigr)+\sum_{k\le Y}\frac{\chi(k)E_0(X/k)}{k^\sigma}
-\sum_{m\le z}\frac{R_\sigma(X/m)}{m^\sigma}.
\end{multline*}
By (i) and (ii), the last two sums are at most $\sum_{k\le Y}k^{-\sigma}(X/k)^{-\sigma}=YX^{-\sigma}$ and
$\sum_{m\le z}m^{-\sigma}q(X/m)^{-\sigma}\le qzX^{-\sigma}=qX^\delta/Y$ in absolute value. With $S\ge1$ and
(iii), this proves the lemma, even with $4$ replaced by $1$.
\end{proof}

The formal proof bounds $|E_0(u)|$ by $4u^{-\sigma}$ for non-integral $u$; the lemma is stated in that form
because the numerical step below is certified for it.

\begin{proposition}\label{prop:L1}
Let $c>0$ and $T\ge2$ with $c/T\le1/2$, and suppose that the zero gap $(c,T)$ holds. Let $A>0$ with $AT\ge e$,
let $V\ge1$ be an integer, and put $\sigma_T=1-c/T$, $P=c(1+\log(2VAT)/T)$ and
\[
\varepsilon=e^c\Bigl(\frac{(AT)^{c/T}V^{c/T}}{cA}+\frac4{(AT)^{\sigma_T}}\Bigr)+\frac{e^P}{cA}
+\frac{4e^P}{V^{\sigma_T}}+\frac{e^P}{AT}.
\]
If $\varepsilon<1$, then every primitive real character $\chi$ of conductor $q\ge3$ with
$\log q\ge T$ satisfies $\sqrt q\,L(1,\chi)\ge(1-\varepsilon)ce^{-P}e^{T/2}/T$.
\end{proposition}

\begin{proof}
Recall $\ell=\log q$. Let $\delta=c/\ell\le1/2$ and $\sigma=1-\delta\ge\sigma_T$. By the zero gap, $L(t,\chi)\ne0$ for
$\sigma<t<1$. Apply \cref{lem:L1} with $Y=\lceil Aq\ell\rceil$ and $X=VY$. Since $Aq\ell\le Y\le2Aq\ell$,
\[
X^\delta\le(2VAq\ell)^{c/\ell}=\exp\Bigl(c+c\,\frac{\log(2VA\ell)}\ell\Bigr)\le e^P,
\]
because $\log(2V)/\ell$ and, as $A\ell\ge e$, $\log(A\ell)/\ell$ decrease in $\ell$. Using $q^\delta=e^c$, $A\ell\ge AT\ge1$, $\sigma\ge\sigma_T$
and that $(A\ell)^{c/\ell}$ and $V^{c/\ell}$ decrease in $\ell$, the four terms of $E$ are bounded by the terms of $\varepsilon$:
$q(Y+1)^{-\sigma}\le e^c(A\ell)^{-\sigma}$, so the first term is at most
$e^c(A\ell)^{\delta}V^\delta/(cA)+4e^c(A\ell)^{-\sigma}$; the second is at most $e^P(\ell/c)/(A\ell)$; the
third equals $4V^{-\sigma}Y^\delta\le4V^{-\sigma_T}e^P$; the fourth is at most $e^P/(A\ell)$. Hence
$E\le\varepsilon$ and $L(1,\chi)\ge\delta X^{-\delta}(1-\varepsilon)\ge(c/\ell)e^{-P}(1-\varepsilon)$. Finally,
$\sqrt q=e^{\ell/2}$ and $e^{\ell/2}/\ell$ increases for $\ell\ge2$.
\end{proof}

Let $w_K\in\{2,4,6\}$ be the number of roots of unity in $K=K_d$. Dirichlet's class number formula
\cite[Ch.~6, (15)]{Davenport}, $2\pi h(-d)=w_K\sqrt d\,L(1,\chi_{-d})$, gives
\begin{equation}\label{eq:cnf}
\pi h(-d)\ge\sqrt d\,L(1,\chi_{-d}).
\end{equation}

\begin{theorem}\label{thm:lean}
Let $N\in\{100,200,500,1000,2000\}$ and let $T_N$ be as in Table~\ref{tab:TN}. Every imaginary quadratic
field $K$ with $\log|d_K|\ge T_N$ has $h_K>N$.
\end{theorem}

\begin{table}[ht]
\caption{The thresholds $T_N$ of \cref{thm:lean} and the zero gaps $(c_N,T_N)$ used in its proof ($e^{T_N}$
and $c_N$ rounded).}
\label{tab:TN}
\begin{tabular}{rlll}
$N$ & $T_N$ & $e^{T_N}$ & $c_N$\\\hline
$100$ & $28.6787$ & $2.8510\cdot10^{12}$ & $0.005506$\\
$200$ & $30.1458$ & $1.2364\cdot10^{13}$ & $0.005622$\\
$500$ & $32.0058$ & $7.9422\cdot10^{13}$ & $0.005759$\\
$1000$ & $33.4757$ & $3.4539\cdot10^{14}$ & $0.005861$\\
$2000$ & $35.0131$ & $1.6069\cdot10^{15}$ & $0.005965$
\end{tabular}
\end{table}

\begin{proof}
Let $K=K_d$ with $\log d\ge T_N$; then $d>8$. Each $T_N$ is a grid point of the proof of \cref{thm:gap}, and
$c_N$ is the minimum of the constants of the intervals contained in $[T_N,\infty)$; hence the zero gap
$(c_N,T_N)$ holds. Take $A=10^6$ and $V=10^4$ in \cref{prop:L1}. Then $\varepsilon<7.8\cdot10^{-4}$,
and $(1-\varepsilon)c_Ne^{-P}e^{T_N/2}/T_N>N\pi$; the margin in the logarithm is between $0.009$ ($N=500$) and
$0.072$ ($N=2000$) (\texttt{computations/lean\_h\_cert.txt}). By \cref{prop:L1} and \eqref{eq:cnf},
$\pi h(-d)\ge\sqrt d\,L(1,\chi_{-d})>N\pi$.
\end{proof}

\begin{remark}\label{rem:halfplane}
Suppose that no Dirichlet $L$-function vanishes in $\operatorname{Re}s>31/32$; each of the claims of
\cite{OAI:The-Quasi-Riemann-Hypothesis-September-30-2026,OAI:The-Quasi-Riemann-Hypothesis-October-5-2026,CuretonThorner2026}
implies this. Then $1-\beta\ge1/32$ for every real zero $\beta$, so the zero gap $(c,T)$ holds for every
$c\le T/32$. Taking $A=V=10^3$ and $c=0.595$ (resp.\ $c=0.652$) in \cref{prop:L1}, we obtain from
\eqref{eq:cnf} that $h(-d)>2000$ for $\log d\ge27.2$ (resp.\ $h(-d)>10^5$ for $\log d\ge35.4$)
(\nolinkurl{computations/qrh_cert.txt}); here $e^{27.2}<6.5\cdot10^{11}$ and $e^{35.4}<2.4\cdot10^{15}$. A
larger $c$ does not help, since $ce^{-P}$ is largest at $c=(1+\log(2VAT)/T)^{-1}$, about $0.6$ here; the
half-planes $\operatorname{Re}s>7/8$ and $\operatorname{Re}s>11/12$ give the same thresholds. For $h\le2000$
the sieve would then be needed only on $[6\cdot2^{28},6.5\cdot10^{11}]$, inside the range $d<2^{40}$ of the
tables of \cite{MosunovJacobson2016}. For $h\le10^5$ the bound is of the size of the range of our sieve for
$h\le2000$; we have not attempted the computation below it. Since the half-plane also excludes complex
zeros, the argument of Littlewood \cite{Littlewood1928} gives $h(-d)\gg\sqrt d/\log\log d$ with an effective
implied constant, as noted in \cite[\S1]{OAI:The-Quasi-Riemann-Hypothesis-October-5-2026}. Lamzouri, Li and
Soundararajan \cite[Thm.~1.5, Cor.~1.3]{LLS2015} made this bound explicit under GRH for Dirichlet
$L$-functions; their proof uses that the zeros lie on the line $\operatorname{Re}s=1/2$. We have not made it
explicit under the half-plane.
\end{remark}

\section{Small discriminants: the Mosunov--Jacobson tables}\label{sec:tables}

Mosunov and Jacobson \cite{MosunovJacobson2016} computed the class groups of all imaginary quadratic fields
with $d<2^{40}$. Their tabulation does not assume the extended Riemann hypothesis for Dirichlet $L$-functions:
for $d\not\equiv7\pmod8$ the class numbers are computed unconditionally from class number formulas involving
theta-series, and for $d\equiv7\pmod8$ the class groups are verified with an identity derived from the
Eichler--Selberg trace formula \cite[\S\S3, 5]{MosunovJacobson2016}. We use the class
numbers for $d<6\cdot2^{28}\approx1.61\cdot10^9$.

\subsection{The data}
{\sloppy The tables are distributed by the LMFDB \cite{LMFDB} as files $\mathtt{cl}r\mathtt{mod}m.k\mathtt{.gz}$ for
$(r,m)\in\{(3,8),(7,8),(4,16),(8,16)\}$; these four classes contain every $d$. The file $k$ for $(r,m)$ covers
the $d\equiv r\pmod m$ with $k\cdot2^{28}\le d<(k+1)2^{28}$. Each line holds a gap, the class number and the
invariant factors of the class group, and $d$ advances by $m$ times the gap from the first $d\equiv r\pmod m$
in the range. We determined this format from the files and use the $24$ files with $k\le5$ (SHA-256 in
\nolinkurl{computations/mj_data.sha256}). The script \nolinkurl{computations/mj_extract.py} lists the pairs
$(d,h(-d))$ with $h(-d)\le2000$: there are $16{,}453{,}493$, all with $d<2^{30}$
(\nolinkurl{computations/mj_extract_h2000.txt}).\par}

\subsection{Checks}
\begin{enumerate}
\item[(a)] \emph{Completeness.} In each file the gaps after the first are positive, every decoded $d$ lies in the range of
the file and $-d$ is a fundamental discriminant, and the number of entries equals the number of fundamental
discriminants $-d$ with $d\equiv r\pmod m$ in that range, counted with a separate squarefree sieve. Hence
each file lists every such $d$ exactly once (\nolinkurl{computations/mj_completeness.py}; deleting one line
of a file makes the check fail).
\item[(b)] \begin{sloppypar}Random entries were compared with PARI/GP \cite{PARI}: $160$ entries of all sizes and
$2000$ with $200<h\le1000$ (\nolinkurl{qsieve/TESTING.md}), and $2000$ with $1000<h\le2000$
(\nolinkurl{computations/random_check_h2000.txt}). The class numbers agree with \texttt{qfbclassno}, which
the PARI/GP documentation states to be unconditionally correct for $d<2\cdot10^{10}$ \cite{PARI}.\end{sloppypar}
\item[(c)] The fields with $h\le100$ in the tables are exactly those found by computing $h(-d)$ with \texttt{qfbclassno}
for all $d\le2.5\cdot10^6$ ($42{,}272$ fields); this computation reproduces the counts and the largest $d$
for each $h\le100$ in Watkins' Table~4 \cite{Watkins2004} (\nolinkurl{computations/watkins_table4.txt}).
\item[(d)] The sieve of \cref{sec:sieve}, followed by the post-processing of the uncertified $d$
(\nolinkurl{qsieve/postproc.gp}), finds every field with $h(-d)\le N$ in a given range. Run with
$N=1000$ on $[2^{24},2^{28})$ and with $N=2000$ on $[2^{28},6\cdot2^{28})$, it returns exactly the
$1{,}490{,}910$ and $537{,}404$ fields of the tables in these ranges (\nolinkurl{qsieve/TESTING.md}).
\end{enumerate}
Further comparisons with \cite{CremonaSutherland2024,HJKMP2019} are in \cref{sec:results}.

By (c) and (d), the fields with $h\le1000$ and $2^{24}\le d<2^{28}$, those with $h\le2000$ and
$2^{28}\le d<6\cdot2^{28}$, and those with $h\le100$ and $d\le2.5\cdot10^6$ have been determined by a
computation independent of the tables. For the other $d<6\cdot2^{28}$ our result rests on the class numbers
of \cite{MosunovJacobson2016}, spot-checked by (b).

\section{The sieve}\label{sec:sieve}

On $6\cdot2^{28}\le d\le1.65\cdot10^{15}$ we show $h(-d)>2000$ by the sieve of Watkins
\cite[\S8]{Watkins2004}, which bounds $h(-d)$ from below by counting reduced forms with small first
coefficient. Throughout this section $N$ is the bound to be certified (in the runs, $N=2000$), $-d$ is a
fundamental discriminant with $d>4$, $\chi=\chi_{-d}$, and $X=\sqrt{d/4}$. A prime $p$ is \emph{split}, \emph{ramified} or \emph{inert} according as $\chi(p)=1$, $0$ or
$-1$.

\subsection{The lower bound}
For $n\ge1$ let $m_d(n)$ be the number of reduced forms of discriminant $-d$ with first coefficient $n$. Since
$-d$ is fundamental, every form of discriminant $-d$ is primitive; every class contains exactly one reduced
form, so $h(-d)=\sum_nm_d(n)$.

\begin{lemma}[{\cite[(2), Lemma~3]{Watkins2004}}]\label{lem:minima}
For $n\le X$, $m_d(n)=\prod_{p^k\parallel n}e(p^k)$, where $e(p^k)=2$ if $p$ splits, $e(p)=1$ if $p$
ramifies, and $e(p^k)=0$ otherwise.
\end{lemma}

\begin{proof}
Let $-n<b\le n$ with $b^2\equiv-d\pmod{4n}$. The form with coefficients $n$, $b$ and $(b^2+d)/4n$ has
discriminant $-d$, and $(b^2+d)/4n\ge d/4n\ge n$, with equality only if $b=0$ and $d=4n^2$, which is
impossible for $d>4$; so it is reduced. Conversely, a reduced form with first coefficient $n$ has middle
coefficient $b$ with $-n<b\le n$. Hence $m_d(n)$ is the number of $b\bmod2n$ with $b^2\equiv-d\pmod{4n}$. This number is
multiplicative in $n$ by the Chinese remainder theorem. For an odd prime power $p^k$ it is $1+\chi(p)$ if
$p\nmid d$, and, since $-d$ is fundamental, $1$ if $p\mid d$ and $k=1$, and $0$ if $p\mid d$ and $k\ge2$. For
$n=2^kn'$ with $n'$ odd, the factor at $2$ is the number of $b\bmod2^{k+1}$ with $b^2\equiv-d\pmod{2^{k+2}}$;
it is $1+\chi(2)$ if $d\equiv3\pmod4$, since odd squares are $\equiv1\pmod8$ and an odd square modulo
$2^{k+2}$ has exactly $4$ square roots there, that is, $2$ modulo $2^{k+1}$; and if $d=4m$ with
$m\equiv1,2\pmod4$, putting $b=2c$ gives $c^2\equiv-m\pmod{2^k}$, which has one solution $c\bmod2$ for $k=1$ and
none for $k\ge2$.
\end{proof}

Let $\mathcal P$ be a finite set of primes, each labelled split or ramified, and $Y\le X$. Let
$\mu_{\mathcal P}(Y)$ be the sum of $\prod_{p^k\parallel n}e(p^k)$ over the $n\le Y$ all of whose prime factors
lie in $\mathcal P$ ($n=1$ included), with $e$ as in \cref{lem:minima} taken according to the labels. If the
labels are the true types, then
\begin{equation}\label{eq:mu}
h(-d)\ge\mu_{\mathcal P}(Y).
\end{equation}
The function $\mu_{\mathcal P}(Y)$ is nondecreasing in $\mathcal P$ and in $Y$, and it does not increase if a
prime $p\in\mathcal P$ labelled split is replaced by a prime $p'>p$, $p'\notin\mathcal P$, labelled split: replacing $p'^k$ by $p^k$ maps the terms
for the new set injectively to equal terms for $\mathcal P$.

\subsection{The certificates}
Let $Q=16\cdot9\cdot5\cdot7\cdot11\cdot13\cdot17\cdot19\cdot23\cdot29=155{,}272{,}637{,}520$. The values
$\chi(p)$ for $p\le29$ depend only on $d\bmod Q$, and some classes mod $Q$ contain no fundamental $-d$. Let
$G_1,G_2,\dots,G_{40}$ be the consecutive blocks of $35$ primes starting at $31$, and $v_j=\max G_j$. The
range is divided into segments, the intersections of the range with the intervals $[2^k,2^{k+1})$; for a
segment with least element $a$ put $Y_a=\lfloor\sqrt{a/4}\rfloor\le X$. For each pattern of values $\chi(p)$, $p\le29$, let $\mathcal P_0$ be the set of
split primes $p\le29$ and ramified primes $p\le\min(29,Y_a)$, and let $s$ be the least integer $s\ge0$ with
$\mu_{\mathcal P_0\cup\{v_1,\dots,v_s\}}(Y_a)>N$, the $v_j$ labelled split ($s=\infty$ if there is no
such $s\le40$). A number $d$ is \emph{certified} if
one of the following holds.
\begin{enumerate}
\item[(i)] $s=0$, that is, $\mu_{\mathcal P_0}(Y_a)>N$.
\item[(ii)] $s\le40$, and for each $j\le s$ some prime of $G_j$ splits.
\item[(iii)] For one of the nine primes $b^*$ that are the $8$th, $16$th, $32$nd, \dots, $1024$th and
$1400$th prime after $29$, with $t$ the number of split primes in $(29,b^*]$ and
$W_0=\mu_{\mathcal P_0}(Y_a/(b^*)^2)$, we have $(2t^2+1)W_0>N$.
\item[(iv)] $\mu_{\mathcal P}(X)>N$ for the set $\mathcal P$ of all split or ramified primes $p\le\min(X,2\cdot10^6)$.
\end{enumerate}

\begin{lemma}\label{lem:certified}
If $d$ is certified, then $h(-d)>N$.
\end{lemma}

\begin{proof}
(i) and (iv) are \eqref{eq:mu}. (ii): let $p_j\in G_j$ split. The $p_j$ are distinct, larger than $29$
and $p_j\le v_j$; by the monotonicity of $\mu$, $h(-d)\ge\mu_{\mathcal P_0\cup\{p_1,\dots,p_s\}}(Y_a)\ge
\mu_{\mathcal P_0\cup\{v_1,\dots,v_s\}}(Y_a)>N$. (iii): let $p_1,\dots,p_t$ be the split primes in $(29,b^*]$. For
each $n$ counted in $W_0$, the numbers $n$, $np_i$ and $np_ip_{i'}$ ($i<i'$) are at most $Y_a\le X$; they are
pairwise distinct, also for different $n$, since the prime factors of $n$ are at most $29$. By
\cref{lem:minima}, $m_d(np_i)=2m_d(n)$ and $m_d(np_ip_{i'})=4m_d(n)$, so these numbers contribute
$(1+2t+4\binom t2)\,m_d(n)=(2t^2+1)\,m_d(n)$ to $h(-d)=\sum_nm_d(n)$. Summing over $n$ gives
$h(-d)\ge(2t^2+1)W_0>N$.
\end{proof}

The tests are applied in the order (i)--(iv). The $d$ that pass none of them would be written to a file and
their class numbers computed with PARI/GP (\nolinkurl{qsieve/postproc.gp}); in the runs below there were
none. The values of $d$ for which $-d$ is not fundamental need no treatment; they are skipped where their
class mod $Q$ shows this, and otherwise processed like the others, which is harmless.

\subsection{Implementation}
The program \nolinkurl{qsieve/qsieve.c} computes the patterns, the numbers $s$ and the $W_0$ for each
segment. If a segment is at least $Q$ long, it runs over all classes $r\bmod Q$ and over the $d\equiv r$ in
the segment; otherwise over all $d$ in it. The work is split into units (a segment and a class modulo
$23\cdot29$, or a subinterval) distributed over processes, and each process logs the units it has
finished. For classes with at least $10$ values of $d$, test (ii) is evaluated for all of them at once on
bit vectors (the periodic patterns of $\chi(p)$, $p\in G_j$, along $d\equiv r\pmod Q$), with AVX-512
instructions.

The program was checked as follows (\nolinkurl{qsieve/TESTING.md}):
\begin{enumerate}
\item[(a)] a self-test compares all tabulated values of $\chi(p)$ with directly computed Kronecker symbols
($666{,}511$ random $d<2^{62}$), and the three implementations of the bit-vector test with the test applied
to each $d$ ($9{,}199{,}491$ cases);
\item[(b)] with every class examined individually, the numbers of $d$ and of classes visited agree with
their exact values;
\item[(c)] on $3\le d\le2.5\cdot10^6$ with $N=100$ it returns the $42{,}272$ fields of \cref{sec:tables}(c),
and it reproduces the tables on $[2^{24},2^{28})$ with $N=1000$ and on $[2^{28},6\cdot2^{28})$ with $N=2000$
(\cref{sec:tables}(d));
\item[(d)] in test runs, every certification by (iii) was confirmed by (iv), and the counters of the runs
agree between the implementations of the bit-vector test.
\end{enumerate}

\subsection{The runs}
The runs for $N=2000$ were made on an AMD Ryzen Threadripper PRO 7995WX (96 cores) with $188$ processes:
\begin{center}\small
\begin{tabular}{lrrrr}
range of $d$ & units & time (s) & $d$ examined & remaining\\\hline
$[6\cdot2^{28},2^{40}]$ & $13{,}367$ & $3{,}778$ & $2.20\cdot10^{11}$ & $0$\\
$[2^{40},1.65\cdot10^{15}]$ & $7{,}337$ & $10{,}142$ & $1.55\cdot10^{14}$ & $0$
\end{tabular}
\end{center}
Here ``$d$ examined'' counts the $d$ handled individually or in bit vectors; the others lie in classes
certified by (i) or contain no fundamental $-d$. All processes finished all their units, and no $d$ remained uncertified. Logs are in
\nolinkurl{qsieve/logs/gap2000} and \nolinkurl{qsieve/logs/main2000}. By \cref{lem:certified}:

\begin{proposition}[computation]\label{prop:sieve}
No imaginary quadratic field with $6\cdot2^{28}\le d\le1.65\cdot10^{15}$ has class number at most $2000$.
\end{proposition}

Since $e^{T_N}<1.65\cdot10^{15}$ for all $N$ in \cref{tab:TN} ($e^{T_{2000}}\approx1.6069\cdot10^{15}$), the
sieved range meets the range of \cref{thm:lean} for each $N$ (\nolinkurl{computations/lean_threshold_match.txt}).

\section{Results and comparisons}\label{sec:results}

\begin{proof}[Proof of \cref{thm:main}]
Let $K=K_d$ with $h(-d)\le2000$. By \cref{thm:lean}, $d<e^{T_{2000}}<1.65\cdot10^{15}$, and by
\cref{prop:sieve}, $d<6\cdot2^{28}$. By \cref{sec:tables}, the fields with $d<6\cdot2^{28}$ and
$h(-d)\le2000$ are the $16{,}453{,}493$ fields of the list, and the largest $d$ among them is
$1{,}062{,}350{,}083$, with $h=1978$.
\end{proof}

The same argument with $N<2000$ gives the lists for $h\le N$. \Cref{tab:counts} gives their size; the number
of fields and the largest $d$ for each $h\le2000$ are in \nolinkurl{computations/h2000_table.txt}, and the
list itself is in the data set (\hyperref[sec:data]{Data and code availability}).

\begin{table}[ht]
\caption{Imaginary quadratic fields with $h(-d)\le N$.}\label{tab:counts}
\begin{tabular}{rrrr}
$N$ & fields & largest $d$ & $h(-d)$\\\hline
$100$ & $42{,}272$ & $2{,}383{,}747$ & $98$\\
$200$ & $167{,}157$ & $8{,}377{,}363$ & $198$\\
$500$ & $1{,}028{,}749$ & $55{,}171{,}603$ & $494$\\
$1000$ & $4{,}115{,}897$ & $227{,}932{,}027$ & $996$\\
$2000$ & $16{,}453{,}493$ & $1{,}062{,}350{,}083$ & $1978$
\end{tabular}
\end{table}

\subsection{Comparisons}
{\sloppy
\emph{$h\le100$.} The number of fields and the largest $d$ for each $h\le100$ agree with Watkins' Table~4
\cite{Watkins2004} (\cref{sec:tables}(c)), with the sequences A046125 \cite{OEIS-A046125} (number of fields)
and A038552 \cite{OEIS-A038552} (largest squarefree $k$ with $h(\Q(\sqrt{-k}))=h$), and the $42{,}271$ pairs
$(d,h(-d))$ with $d\le2{,}250{,}000$ agree with the number fields in the LMFDB \cite{LMFDB}, which are
complete in that range (\nolinkurl{computations/external_compare_h1000.txt},
\nolinkurl{computations/lmfdb_check.txt}).

\emph{Cremona and Sutherland.} In this paragraph $D<0$ is the discriminant of an imaginary
quadratic order, not necessarily fundamental: $D=f^2D_0$ with $D_0$ fundamental and $f\ge1$, $\mathcal O_D$ is
the order of discriminant $D$ in $K=\Q(\sqrt{D_0})$, $\mathcal O_K$ the ring of integers of $K$, and $h(D)$ the
class number of $\mathcal O_D$. Under GRH, Cremona and Sutherland \cite[\S5]{CremonaSutherland2024} list the
$6{,}450{,}424$ such $D$ with $h(D)\le1000$ (file \nolinkurl{cmdiscs1000.txt} of their repository, commit
\texttt{1cbb56fc}). We have
\[
h(D)=\frac{h(D_0)f}{[\mathcal O_K^\times:\mathcal O_D^\times]}\prod_{p\mid f}\Bigl(1-\Bigl(\frac{D_0}p\Bigr)\frac1p\Bigr)
\]
\cite[(2.1)]{CremonaSutherland2024}; in particular $h(D_0)\mid h(D)$. From our fields with $h\le1000$ this
formula generates every $D$ with $h(D)\le1000$, and the resulting pairs $(D,h(D))$ coincide with theirs
(\nolinkurl{computations/compare_cs1000.txt}). As \cite[\S5]{CremonaSutherland2024} explain, their list can
also be obtained from the same tables \cite{MosunovJacobson2016}, below a bound for $|D_0|$ that holds under GRH \cite[Prop.~19]{CremonaSutherland2024}; it comes from a lower
bound for $L(1,\chi_{D_0})$ that assumes GRH for Dirichlet $L$-functions \cite[Cor.~1]{LanguascoTrudgian2022}.
The agreement checks our extraction; with \cref{thm:main} it shows that their list is complete
unconditionally.

\emph{Holmin, Jones, Kurlberg, McLeman and Petersen.} Under GRH, \cite{HJKMP2019} give the number of
imaginary quadratic fields with class number $h$ for every odd $h<10^6$. For the $1000$ odd $h\le1999$ these
numbers agree with ours ($1{,}712{,}058$ fields in total; \nolinkurl{computations/compare_hk2000.txt}), so
they hold unconditionally. Their computation assumes GRH in two places \cite[\S6]{HJKMP2019}: for Dirichlet
$L$-functions, in the bound for $d$ and in the approximation of $h(-d)$; and, through the class groups computed
with PARI/GP, for the $L$-functions of the nontrivial characters of the class group of $K_d$, which are Hecke
$L$-functions and in general not Dirichlet $L$-functions. The latter hypothesis is needed to know that prime
ideals of small norm generate the class group \cite[Thm.~2.4]{GrenieMolteni2018}.

\emph{Garcia, Greni\'e, Lee and Molteni.} Assuming GRH for the Dedekind zeta function of $K_d$,
\cite[(10)]{GGLM2026} show that $h(-d)\le2000$ implies $d\le8{,}213{,}527{,}466$, and they determine the fields
below this bound with PARI/GP. The list they publish \cite[p.~4]{GGLM2026} (file
\nolinkurl{classes-list-2000.txt}, downloaded on 22 October 2026) consists of the same $16{,}453{,}493$ pairs
$(d,h(-d))$ as ours, so it holds unconditionally
(\nolinkurl{computations/compare_gglm2000.txt}).

\emph{Platt and Trudgian.} Counting the reduced forms of every fundamental discriminant $-d$ with
$d\le10^{10}$ \cite[\S3.1]{PlattTrudgian2026}, Platt and Trudgian find $16{,}453{,}493$ fields with
$h(-d)\le2000$, the largest $d$ being $1{,}062{,}350{,}083$ \cite[\S4.2]{PlattTrudgian2026}. Their
computation does not use \cite{MosunovJacobson2016}, so this checks the number of fields and the largest $d$
in \cref{tab:counts} independently.

\emph{Jaswal.} Assuming the second statement of \cref{thm:gap}, Jaswal \cite{Jaswal2026} determines the
fields with class number at most $3000$; for $d<1.98\cdot10^{10}$ his computation is unconditional and does
not use \cite{MosunovJacobson2016}. His fields with $h(-d)\le2000$ are those of our list
\cite[App.~A.3]{Jaswal2026}.\par}

\section{Formal verification}\label{sec:lean}

\Cref{thm:gap,thm:lean} are proved in Lean~4 \cite{deMouraUllrich2021} with Mathlib
\cite{Mathlib2019}. The files are in the directory \nolinkurl{lean/} of the code repository (\hyperref[sec:data]{Data and code availability}).

\subsection{The statements}
For $N=2000$ the Lean theorem reads
\begin{flushleft}\footnotesize\ttfamily
theorem B4.ClassNumber.classNumber\_gt\_2000 (K : Type*) [Field K]\\
\quad[NumberField K] [IsTotallyComplex K] (hK : finrank $\mathbb Q$ K = 2)\\
\quad(h : (350131/10000 : $\mathbb R$) $\le$\\
\qquad Real.log |(NumberField.discr K : $\mathbb R$)|) :\\
\quad 2000 < NumberField.classNumber K
\end{flushleft}
and likewise for $N=100,200,500,1000$ with the thresholds of \cref{tab:TN}. The hypotheses say that $K$ is a
totally complex number field of degree~$2$, that is, an imaginary quadratic field, and every notion in the
statement is a Mathlib definition. This is \cref{thm:lean}. \Cref{thm:gap} is
\texttt{B4.Chi8.zero\_gap\_all} (the bound $5811\cdot10^{-7}$) and \texttt{B4.Chi8.zero\_gap\_large};
there $\chi$ is a Dirichlet character in Mathlib's sense with real values, primitive and nontrivial, and
$L(s,\chi)$ is Mathlib's \texttt{LFunction}. The zero gap $(c,T)$ of the Notation is the predicate
\texttt{B4.Chi8.GapFromAll c T}. The certificates are proved for the predicate \texttt{B4.GapFrom c T}, which
assumes in addition $q\ne8$; \cref{lem:chi8}, \texttt{B4.Chi8.re\_LFunction\_pos}, removes this hypothesis.
It is proved there in the form $\operatorname{Re}L(\sigma,\chi)>0$ for all $\sigma>0$, for every nontrivial
character modulo $8$.

\subsection{The proof}
The proof consists of $20$ files, about $7{,}400$ lines, that follow
\cref{sec:zerogap,sec:L1}:
\begin{center}\footnotesize
\begin{tabular}{ll}
\texttt{B4ExplicitBias}, \texttt{B4PrimeMass}, \texttt{B4Analytic} & \cref{lem:bias,lem:primes}\\
\texttt{B4MasterAlgebra}, \texttt{B4ZeroGap}, \texttt{B4ElementaryGap} & \cref{prop:master}\\
\texttt{B4Interval}, \texttt{B4Chain}, \texttt{B4Certificates} & \cref{thm:gap}\\
\texttt{B4Chi8} & \cref{lem:chi8}, \cref{thm:gap} for $q=8$\\
\texttt{B4LSeriesCont}, \texttt{B4ZetaEM}, \texttt{B4Sign}, & \\
\quad\texttt{B4Hyperbola}, \texttt{B4L1} & \cref{lem:L1}\\
\texttt{B4ClassBound} & \cref{prop:L1} for $q\ne8$\\
\texttt{B4ClassCert} & numerical step\\
\texttt{B4QuadSplit}, \texttt{B4QuadZeta}, \texttt{B4ClassNumber} & \eqref{eq:cnf}, \cref{thm:lean}
\end{tabular}
\end{center}

\emph{The upstream proof.} The inequalities \eqref{eq:upper}--\eqref{eq:pivot}, the formula used in
\cref{lem:bias} for $-L'(s,\chi)/L(s,\chi)$, and the bound $\sum_{m\le n}\Lambda(m)/m\ge\log n-1$ are taken
unmodified from the Lean formalization of \SZ{} in \cite{OpenAIMath2026}, directory
\nolinkurl{lean/OAI/NumberTheory/SiegelZeros} at commit \texttt{301488868}. We do not use
its final theorem, which gives Theorem~1 of \SZ{} with an unspecified constant.

\emph{Where the formalization departs from \SZ{}.} The interpolation step is not formalized as
\SZ[Lemma~3], neither its statement for
general degrees nor its proof (a supporting face of a lattice polytope and Lagrange interpolation). What is
formalized is the conclusion of
\SZ[Cor.~4] for our box: the rows $R_\alpha$ of \cref{sec:det}, with $B$, $H$, $U$ as there,
with $\alpha_1\le32H^{2/3}U$ and $\alpha_2,\alpha_3\le32H^{-1/3}U$ span
(\texttt{actual\_biquadratic\_rectangle\_span}). Its proof passes to the dual statement, that a nonzero
polynomial supported on $\{0,\dots,B-1\}^4$ cannot vanish at $(1,1,1,1)$ to all these orders along
the three invariant derivations, and proves it by a multiplicity estimate
(\texttt{uniform\_rectangular\_multiplicity}). This conclusion enters only the pivot bounds
\eqref{eq:pivot} (\texttt{global\_source\_rectangle\_pivot\_bounds}). That $M=B^4$ rows are retained and
their determinant $\Delta$ is nonzero is proved without it: the rows with $\alpha_2=\alpha_3=0$ and $\alpha_1<M$ already span, by
the Vandermonde determinant of the $M$ distinct numbers $\vartheta_n$ (\texttt{jet\_weight\_ball\_span\_eq\_top}).

\emph{Other libraries.} Mathlib provides $\theta(x)\le x\log4$, the convexity of $\log\Gamma$, the
continuation of $L(s,\chi)$, and the residue of the Dedekind zeta function $\zeta_K(s)$ at $s=1$ in terms of
$h_K$, the regulator and the number of roots of unity. PrimeNumberTheoremAnd \cite{PNTplus} provides the
integral representation of $\zeta(\sigma)$ for $0<\sigma<1$ and the bound for the prime powers in
\cref{lem:primes}(c). TauCeti \cite{TauCeti} provides the Kronecker character of a quadratic field and the
splitting of primes in it, from which we prove $\zeta_K(s)=\zeta(s)L(s,\chi_D)$ for $\operatorname{Re}s>1$,
using the Euler products of AINTLIB \cite{AINTLIB}; together with the residue formula this gives the class
number formula. Lean~4.34.1 and the dependency versions are those pinned by the upstream project (Mathlib
commit \texttt{d13f23b}). The code repository is a Lake project that requires Mathlib, LeanCert and
LeanArchitect at these commits. The files of the upstream formalization, of PrimeNumberTheoremAnd, of
TauCeti and of AINTLIB that the proof imports are copied into \nolinkurl{lean/vendor/};
\nolinkurl{lean/vendor/check.sh} compares them byte for byte with the upstream commits.

\emph{Numerical steps.} The $314$ inequalities in the proof of \cref{thm:gap} and the inequalities
$(1-\varepsilon)c_Ne^{-P}e^{T_N/2}/T_N>N\pi$ in the proof of \cref{thm:lean}, with $\varepsilon$ and $P$ as in
\cref{prop:L1}, are proved by the interval
arithmetic of LeanCert \cite{LeanCert} in the mode \texttt{trust := kernel}, in which the Lean kernel checks
the evaluation. The Lean files with these inequalities are generated by the scripts
\nolinkurl{computations/lean_grid.gp} and \nolinkurl{computations/lean_h_cert.gp}; the generator plays no role in their correctness.

\subsection{Checks}
\begin{enumerate}\sloppy
\item[(a)] Our files contain no \texttt{sorry}, no \texttt{axiom} and no \texttt{native\_decide}.
\texttt{\#print axioms} reports only \texttt{propext}, \texttt{Classical.choice} and \texttt{Quot.sound}
for the five theorems \texttt{classNumber\_gt\_N} and for \texttt{zero\_gap\_all},
\texttt{zero\_gap\_large} and \texttt{re\_LFunction\_pos} (\nolinkurl{lean/B4Build.txt}).
\item[(b)] \emph{Kernel replay.} The dependency closure of these theorems, $123{,}010$ constants
including the parts of Lean's core library, Mathlib and the other libraries that they use, was checked again
by the Lean kernel, starting from an empty environment (\nolinkurl{lean/B4Replay.lean}, log
\nolinkurl{lean/B4Replay.txt}). The only axioms in the closure are the three above, and the replayed statements are
unchanged.
\item[(c)] \emph{Comparator.} The challenge files \nolinkurl{lean/Challenge/B4ClassNumber.lean} and
\nolinkurl{lean/Challenge/B4ZeroGap.lean} import only Mathlib and state the five theorems
\texttt{classNumber\_gt\_N} and the three theorems above with \texttt{sorry}. For each of them, Comparator
\cite{Comparator} built our proof
in a sandbox, checked that its statements coincide with those of the challenge and that it uses no other
axioms, and had it checked by the Lean kernel and by nanoda \cite{nanoda}, an independent implementation of
the kernel in Rust, which reads the declarations exported by lean4export \cite{lean4export}. Both accepted it. With $T_{1000}$ changed to $33$, or the constant $5811\cdot10^{-7}$ to
$5812\cdot10^{-7}$, in the challenge, Comparator rejected the proof
(\nolinkurl{lean/comparator/B4Comparator.txt}). Comparator reused build products compiled earlier, outside the sandbox, in the same workspace;
Comparator's guarantee against a deliberately adversarial proof assumes a clean build, but the kernel
checks above are unaffected.
\end{enumerate}

\subsection{What the formal proof covers}
{\sloppy The trusted base of \cref{thm:gap,thm:lean} consists of the Lean kernel (with nanoda as a second
implementation, which reads the proof through lean4export), the three axioms, and the Mathlib
definitions occurring in the statements; both kernels include arithmetic on natural-number literals
(GMP). The
upstream formalization, the other libraries and the generated certificates enter only through proofs checked
by the kernel. Three lemma calls in TauCeti were adapted to changed names and arguments in the pinned Mathlib
(\nolinkurl{lean/vendor/patches/TauCeti-pinned-mathlib.patch}); the files of PrimeNumberTheoremAnd carry the
compatibility patch of the upstream project, and the upstream formalization, AINTLIB, Mathlib, LeanCert and
LeanArchitect are unpatched.\par}

In particular, the formal proof does not depend on the text of \SZ{}. It is not an independent proof of
\SZ{}: it follows the argument of \SZ{} except for the interpolation step, and its core,
\eqref{eq:upper}--\eqref{eq:pivot}, is the upstream formalization of that argument, with the interpolation
step proved by a multiplicity estimate as described above. Not formalized are the computations of \cref{sec:tables,sec:sieve}.

\section{Contributions of AI}\label{sec:ai}

\emph{The OpenAI model.} The manuscript \SZ{} and its Lean formalization were produced by an internal model of
OpenAI and released by OpenAI \cite{OAI:Uniform-exclusion-of-Landau-Siegel-zeros-October-1-2026,
OpenAIMath2026}. We took no part in them. Our formal proof uses parts of that formalization unmodified
(\cref{sec:lean}).

\emph{Claude.} The following was done with Claude Opus 5.5 (Anthropic), used through Claude Code, under the
author's direction:
\begin{enumerate}
\item[(a)] the explicit form of the argument of \SZ{} and the lower bound for $L(1,\chi)$
(\cref{sec:zerogap,sec:L1});
\item[(b)] the Lean files of \cref{sec:lean}, the generation of the numerical certificates, the kernel
replay and the Comparator setup;
\item[(c)] the sieve program of \cref{sec:sieve}, the scripts in \nolinkurl{computations/}, and the checks
and comparisons of \cref{sec:tables,sec:results};
\item[(d)] the text of \cref{sec:zerogap,sec:L1,sec:lean,sec:tables,sec:sieve,sec:results} and of this
section, written from the author's outline; the introduction is based on the author's draft.
\end{enumerate}

\emph{The author} posed the problem and made the decisions that shaped the work: to make the constant of \SZ{}
explicit, to base the result on a formal proof rather than on the text of \SZ{}, to extend the computation to
$h\le2000$, and which checks to make, among them the trusted base of the formal proof (the permitted axioms and
the independent kernel nanoda) and the data for the comparisons \cite{CremonaSutherland2024,HJKMP2019}. The
author worked through the proofs in the paper and checked the cited results against their sources. The
author is responsible for the paper.

\emph{What does not depend on trusting the AI.}
\begin{enumerate}
\item[(a)] \Cref{thm:gap,thm:lean} are checked by the Lean kernel and by nanoda (\cref{sec:lean}). What has
to be read by a person is their statements, which use only Mathlib definitions. The author has checked that
they express \cref{thm:gap,thm:lean}.
\item[(b)] The correctness of the sieve rests on \cref{lem:minima,lem:certified}, proved in the paper, and on
the program implementing them. The program was tested against independent data (\cref{sec:tables}(c), (d),
\cref{sec:sieve}); it was not formally verified.
\item[(c)] The list agrees with the published data of \cref{sec:results}, computed by others; the comparison
with \cite{CremonaSutherland2024} checks our extraction from the tables, not the tables themselves.
\item[(d)] The proofs written out in \cref{sec:zerogap,sec:L1,sec:sieve} can be checked by the reader.
\end{enumerate}
The remaining assumptions are the correctness of the tables of \cite{MosunovJacobson2016} where they were
not recomputed (\cref{sec:tables}), of PARI/GP, of the sieve program, and of our extraction script
\nolinkurl{computations/mj_extract.py} (checked by \cref{sec:tables}(a)--(d) and \cref{sec:results}).

\section*{Data and code availability}\label{sec:data}

\emph{Lists.} The lists of the $4{,}115{,}897$ fields with $h\le1000$ and of the $16{,}453{,}493$ fields with
$h\le2000$, one line ``$d$ $h(-d)$'' per field in increasing order of $d$, are the files
\nolinkurl{lists/imq_h1000.txt.xz} and \nolinkurl{lists/imq_h2000.txt.xz} of the code repository below,
released under CC0. The SHA-256 checksums of the uncompressed files are
\begin{center}\scriptsize
\begin{tabular}{ll}
$h\le1000$ & \texttt{075a201c8e7cce6c0b6013423a5d96c6308452afab2fa79882388fc46d8b0c45}\\
$h\le2000$ & \texttt{07100ea3589ab2b490122fa69b438cd5393bdcd2e547dfc29bb3b5c86f9f24c0}
\end{tabular}
\end{center}

\emph{Code.} The Lean files (\cref{sec:lean}), the sieve program and its logs (\cref{sec:sieve}), and the
scripts and outputs cited in the text are in the code repository \url{https://github.com/iwaokimura/cln-imq-2k}, tag \texttt{v1.4} (commit \texttt{e66ab86}). File names in
the text refer to this version.

\emph{Third-party data.} The tables of \cite{MosunovJacobson2016} are available from the LMFDB
\cite{LMFDB}, and \nolinkurl{computations/mj_data.sha256} lists the SHA-256 checksums of the $24$ files we
used.
The data of \cite{CremonaSutherland2024,HJKMP2019,OEIS-A046125,OEIS-A038552} are available from their
authors; our scripts download them or read them from a local copy, and record their checksums. We do not
redistribute third-party data.

\bibliographystyle{amsalpha}
\bibliography{refs}
\end{document}